%% file: GauteroHeusenerAGT.tex
\documentclass[a4paper,leqno,12pt]{amsart}
\usepackage{amssymb,amscd}
\usepackage[latin1]{inputenc}
 \divide\oddsidemargin by 2
\usepackage{enumerate}
\makeatletter
\renewcommand\theenumi{\@alph\c@enumi}
\renewcommand\theenumii{\@alph\c@enumii}
\renewcommand\theenumiii{\@alph\c@enumiii}
\renewcommand\theenumiv{\@alph\c@enumiv}

\let\@listlla\list

\def\list#1#2{\@listlla{#1}{#2\itemsep=2pt\parsep=0pt\topsep=3pt plus 1pt minus 1 pt}}
\newcommand{\tridiagram}[6]{{\par\par \centering
\@picture(120,120)(0,0) \put(30,95){\makebox(0,0)[r]{$#1$}}
\put(90,30){\makebox(0,0)[tl]{$#3$}}
\put(90,95){\makebox(0,0)[l]{$#2$}}
\put(60,102){\makebox(0,0)[b]{$#4$}}
\put(102,60){\makebox(0,0)[l]{$#6$}}
\put(50,50){\makebox(0,0)[tr]{$#5$}} \thinlines
\put(40,95){\vector(1,0){40}} \put(95,80){\vector(0,-1){40}}
\put(25,80){\vector(1,-1){55}}
\endpicture\par\par}\noindent\ignorespaces}
\def\@map#1#2[#3]{\mbox{$#1 \colon #2 \longrightarrow #3$}}
\def\map#1#2{\@ifnextchar [{\@map{#1}{#2}}{\@map{#1}{#2}[#2]}}
\newcommand{\RR}{{\bf R}}

\newcommand{\NM}{{\mathbb N}}

\newcommand{\Aut}[1]{\mbox{\rm Aut}{(#1)}}

\newcommand{\Ker}{\mbox{\rm Ker}}

\newcommand{\mr}{\mathbb{R}}

\newcommand{\mz}{\mathbb{Z}}

\newcommand{\mn}{\mathbb{N}}

\newcommand{\F}[1]{\ensuremath{\mathbb{F}_{#1}}}

\makeatother
\newtheorem{theorem}{Theorem}[section]
\newtheorem{theo}{Theorem}[]
\newtheorem{proposition}[theorem]{Proposition}

\newtheorem{lemma}[theorem]{Lemma}
\theoremstyle{definition}
\newtheorem{definition}[theorem]{Definition}

\newtheorem{remark}[theorem]{Remark}
\newtheorem{convention}[theorem]{Assumption}
\newtheorem{claim}{Claim}

\usepackage{color}

\def\co{\colon\thinspace}
\def\Def{ := }

\usepackage[pdftex]{graphicx,epsfig}

\title[Cohom. carac. of Rel. Hyp., Comb. Thm.]{Cohomological characterization of Relative
Hyperbolicity and Combination Theorem.}
\author{Fran\c{c}ois Gautero, Michael Heusener}
\address{Universit\'e Blaise Pascal,
Campus des C\'ezeaux, Laboratoire de Math\'ematiques, 63177
Aubi\`ere, France} \email{gautero@math.univ-bpclermont.fr,
heusener@math.univ-bpclermont.fr} 


\begin{document}
\maketitle

\begin{abstract}
We give a cohomological characterization of Gromov relative
hyperbolicity. As an application we prove a converse to the
combination theorem for graphs of relatively hyperbolic groups given
in \cite{GauteroDernier}. We build upon, and follow the ideas of,
the work of S. Gersten \cite{Gersten} about the same topics in the
classical Gromov hyperbolic setting.
\end{abstract}

\section*{Introduction}

The celebrated Gromov hyperbolic groups \cite{Gromov} form a central
class of groups in Geometric Group Theory. The paper \cite{Gersten}
gives a cohomological characterization of the Gromov hyperbolicity,
which was, up to then, a purely geometrical notion. As an
application, the author proved the converse to the Bestvina-Feighn
combination theorem \cite{Bestvina-Feighn} for graphs of hyperbolic
groups. More precisely, given a finite graph $\mathcal G$ of
hyperbolic groups and assuming the quasi convexity of the edge
groups in the vertex groups, he proved the necessity of the
so-called ``Annuli Flare'' property for the fundamental group of
$\mathcal G$ being hyperbolic.

Since then, relative hyperbolicity has appeared in Geometric Group
Theory, and is the object of a lot of interest nowadays. Although
already present in \cite{Gromov}, it really grew with Farb's
formulation \cite{Farb}. Among all the definitions which now coexist
\cite{Bowditch, Osin}, two are not equivalent \cite{Sz, bumagin}. We
will term them weak and strong: Gromov relative hyperbolicity is the
strong one \cite{Sz} and, of course, implies weak relative
hyperbolicity \cite{bumagin}. In order to give an illustration of
these two notions, let us just evoke two classical examples:

\begin{itemize}
  \item An example of a strongly relatively hyperbolic group, which is not a
hyperbolic one, is the fundamental group of a hyperbolic, finite
volume manifold with cusps: the relative part
consists of the peripheral subgroups.
  \item An example of a weakly relatively hyperbolic group which is not a strongly
  relatively hyperbolic one is $\mz \oplus \mz$: just put any of the
  $\mz$-factors in the relative part. In the same line of idea, the
  mapping-class groups of compact surfaces are weakly relatively
  hyperbolic \cite{MM}
but not
  strongly relatively hyperbolic in a non-trivial way (as soon as
  the surfaces have sufficiently high complexity) \cite{Aramayona}.
\end{itemize}

We refer the reader to Section \ref{rh} for the definitions about
relative hyperbolicity. A general combination theorem for graphs of
relatively hyperbolic groups, similar to the Bestvina-Feighn theorem
for graphs of hyperbolic groups, has been proven in
\cite{GauteroDernier} (see also \cite{Reeves}). For previous results
in this direction, see \cite{Dahmani,alibegovic}.

The purpose of this short paper is to adapt and extend the above
cited results of \cite{Gersten} to the setting of the strong
relative hyperbolicity. Our first step is to get a
well-suited notion of the $\ell_{\infty}$-cohomology of a group $G$
relative to a family of subgroups $\mathfrak H$. We borrow from
\cite{BE} the definition of the relative cohomology of such a pair
$(G,{\mathfrak H})$. Our first result is stated as follows:

\begin{theo}
\label{sv} The second relative $\ell_\infty$-cohomology of a group
$G$ relative to a family of subgroups $\mathfrak H$, denoted by
$H^2_{(\infty)}(G,{\mathfrak H})$, is well-defined as soon as $G$
admits a finite presentation relative to $\mathfrak H$. If $G$ is
strongly hyperbolic relative to $\mathfrak H$, then
$H^2_{(\infty)}(G,{\mathfrak H})$ strongly vanishes.
\end{theo}

For {\em strong vanishing}, see Definition \ref{gersten}. This
theorem is false in the setting of weak relative hyperbolicity, see \ref{contrexemple}.

Among many other results, the paper \cite{Groves} is also interested
in homological characterizations of the strong relative
hyperbolicity and generalizes some of Gersten's work in
\cite{Gersten}.
The approach
there is however different than our's, in the sense that the authors consider
absolute cycles with compact support (instead of relative cycles
with non-compact support). Another paper dealing with the same subject
is \cite{Asli}, where the authors consider both usual and
bounded cohomology.

The strong vanishing is necessary to get the announced application
about the combination theorem. The application we give below
concerns only semi-direct products of strongly relatively hyperbolic
groups with a free group. The reason is that the subgroups to put in
the relative part are somewhat tedious to describe in the general
case. This would lead to a heavy formulation, without introducing
new interesting phenomena, which all appear in the semi-direct
product case. This semi-direct product case is in some sense a
``generic'' non-acylindrical case, and the most sophisticated one
which might appear as the fundamental group of a graph of groups. If
one wishes to treat semi-direct products with groups which are not
free, one is led to work on $2$-complexes of groups.

The uniform free groups of relatively hyperbolic automorphisms which
appear below were defined in \cite{GauteroDernier}. Definitions are
recalled in Section \ref{CCT}.

\begin{theo}
\label{un cas particulier} Let $G$ be a group which is strongly
hyperbolic relative to a finite collection of subgroups, denoted by
$\mathfrak H$. Let  $\Aut{G, {\mathfrak H}}$ be the group
of relative automorphisms of $(G, {\mathfrak H})$. Let
$\alpha \colon \F{n} \rightarrow \Aut{G, {\mathfrak H}}$ be a
monomorphism from the rank $n$ free group into  $\Aut{G, {\mathfrak
    H}}$. If $G \rtimes_{\alpha} \F{n}$
is finitely generated and strongly hyperbolic relative to a
$\F{n}$-extension of ${\mathfrak H}$ then $\F{n}$ is a uniform free
group of relatively hyperbolic automorphisms of $(G,{\mathfrak H})$.
\end{theo}

The reverse implication is the object of \cite{GauteroDernier}, but
is proved there only in the case of a finitely generated $G$. This
restriction is due to the fact that, in \cite{GauteroDernier},
Farb's approach of relative hyperbolicity is used: this approach
requires the finite generation of $G$. We would like also to point
out that the condition of finite generation that we require here for
$G \rtimes_{{\mathcal A}} \F{n}$ is not necessary if one works with
a slightly restricted kind of $\F{n}$-extensions. Namely, those
which fix (up to conjugacy) each subgroup of $\mathfrak H$. In fact
we feel that this condition of finite generation could be dropped in
the general case, at the expend of some additional
technical work. So, we only give an
indication  (Remark~\ref{moi}) of how one could get rid of this
assumption.

\section{Relative hyperbolicity}
\label{rh}
In \cite{Gersten} the proof of the vanishing theorem for the
$\ell_{\infty}$-cohomology for hyperbolic groups uses in a crucial
way the existence of a linear isoperimetric inequality.
This
designates Osin approach (see \cite{Osin}) of relative hyperbolicity
as the ideal candidate for our purpose.

We recall some basic definitions from \cite[Ch.~2]{Osin}: let $G$ be
a group,  $\mathfrak{H} =(H_\lambda)_{\lambda \in\Lambda}$ a family
of subgroups of $G$ and $X\subset G$. We say that $X$ is a
\emph{relative generating set of $G$ with respect to $\mathfrak H$}
if $G$ is generated by $\big(\cup_{\lambda\in\Lambda}
H_\lambda\big)\cup X$. In the sequel we will always assume that $X$
is symmetric. In this situation $G$ is a quotient of the free
product
\[ F = \big(\ast_{\lambda\in\Lambda} \widetilde H_\lambda\big)\ast F(X)\]
where the groups $ \widetilde H_\lambda$ are isomorphic copies
of
$H_\lambda$ and $F(X)$ is the free group with the basis $X$. Let
us denote by $\mathcal H$ the disjoint union
\[ \mathcal H = \bigsqcup_{\lambda \in
\Lambda} \widetilde H_\lambda \smallsetminus\{1\}\] and by
$(\mathcal H \cup X)^\ast$ the free monoid generated by $\mathcal H
\cup X$. For every $\lambda\in\Lambda$, we denote by $\mathcal
S_\lambda$ the set of all words over the alphabet $\widetilde
H_\lambda \smallsetminus\{1\}$ that represent the identity in
$F$.
The isomorphism $\widetilde H_\lambda\to H_\lambda$ and the
identity map on $X$ can be uniquely extended to a surjective
homomorphism $\epsilon\co F\to G$. We say that $G$ has the
\emph{relative presentation}
\begin{equation}\label{eq:pres}
\langle X,\mathcal H \mid S=1,\; S\in
\mathcal S
=\sqcup_{\lambda\in\Lambda} \mathcal S_\lambda,\; R=1,\;
R\in\mathcal R\rangle
\end{equation}
with respect to $\mathfrak H$, where $\mathcal R \subset (\mathcal H
\cup X)^\ast$, if $\Ker( \epsilon)$ is the normal closure of
$\mathcal R$ in $F$. In the sequel we will write $G=\langle X,
H_\lambda, \lambda\in\Lambda \mid R=1, R\in {\mathcal R}\rangle$ or
$G=\langle X,{\mathcal H} \mid {\mathcal S}, {\mathcal R}\rangle$
for short.

The relative presentation (\ref{eq:pres}) is called \emph{finite}
if
both sets $X$ and $\mathcal R$ are finite. We say that $G$ is  {\em
finitely presented relative to $\mathfrak H$} if there is a finite
relative presentation of $G$ with respect to $\mathfrak H$.

We denote by $l_{{\mathcal H}\cup X}(\gamma)$ the {\em $\mathfrak
H$-relative length} of $\gamma \in G$, defined as the word-length of
$\gamma$ with respect to the system of generators $X \cup {\mathcal
H}$ (that is the minimal number of elements in $X \cup {\mathcal H}$
needed to write $\gamma \in G$). This is nothing else than the
length of a geodesic from the identity to the $\gamma$ in the Cayley
graph $\Gamma(G,X \cup {\mathcal H})$ of $G$ with respect to
$X \cup
{\mathcal H}$.


Let $c$ be an edge-circuit  in $\Gamma(G,X \cup {\mathcal H})$.
Consider a {\em filling} $\Delta$ for $c$, with respect to a finite
relative presentation (\ref{eq:pres}). That is $\partial \Delta$
corresponds to $c$ and decomposes in subcells whose boundary either
corresponds to a relator from $\mathcal S_\lambda$ for some
$\lambda\in\Lambda$ or to a relator from $\mathcal R$. The {\em
$\mathcal R$-relative area of $\Delta$} is the number of $\mathcal
R$-cells that it contains. We denote by $\mathit{Area}^{{\mathcal
R}}(c)$ the {\em $\mathcal R$-relative area of $c$} which is the
minimal $\mathcal R$-relative area of a filling with boundary $c$.

Two functions $f,g \colon \mn \rightarrow \mn$ are
\emph{asymptotically equivalent} if there exist constants $C,K,L$,
$C_0,K_0,L_0$ such that $f(n) \leq g(Cn+K)+Ln$ and $g(n) \leq
f(C_0n+K_0) + L_0 n$.

\begin{definition} \cite{Osin}
Let $G$ be a group which admits a finite relative presentation
$G=\langle X,{\mathcal H} \mid {\mathcal S}, {\mathcal
R}\rangle$
with respect to $\mathfrak H$.

A {\em relative isoperimetric function} for this presentation is a
function $f \co \mn \rightarrow \mn$ such that, for any $n \in \mn$,
for any circuit $c \in \Gamma(G,X \cup {\mathcal H})$ with length
less or equal to $n$, $\mathit{Area}^{{\mathcal R}}(c) \leq f(n)$.

The {\em relative Dehn function for the presentation
$G=\langle X,{\mathcal H} \mid {\mathcal S}, {\mathcal R}\rangle$}
is the smallest relative isoperimetric function for this presentation.

The group $G$ is {\em strongly hyperbolic relative to $\mathfrak H$}
if $G$ is finitely presented with respect to $\mathfrak H$ and the
relative Dehn-function of $G$ for a finite presentation is linear.
\end{definition}

\begin{remark}
Note that two relative Dehn-functions (as defined above) of two
finite relative presentations of $G$ with respect to $\mathfrak H$
are asymptotically equivalent.
Note also that in general not every presentation (\ref{eq:pres}) admits a finite relative Dehn function.
\end{remark}

\section{Relative $\ell_\infty$-cohomology}

Let $(K,L)$ be a CW-pair i.e.\ $K$ is a CW-complex and $L\subset K$
is a subcomplex. We denote by $K^k$ the $k$-skeleton of $K$. Note
that $L^k=K^k\cap L$.

The cellular chain complex $(C_\ast(K,L),\partial)$ is defined in
the usual way (see \cite{Turaev}).
Here and in the sequel we will
work with real coefficients. The chain group $C_k(K,L)= C_k(K) /
C_k(L)$ is usually thought of as the vector space generated by the
$k$-cells of $K\smallsetminus L$.  We denote by
$C^k_{(\infty)}({K},{L})\subset C^k(X,Y)$ the subgroup of bounded
relative $k$-cochains. These cochains correspond exactly to the
$k$-cochains $h\co C_k(K)\to \mr$ which vanish on the $k$-cells of
$L$ and which admit an uniform upper-bound $M_h$ over all $k$-cells
$e^k$ from $K^k$, i.e.\ $h(e^k)=0$ for all $e^k\in L^k$ and
$|h(e^k)| \leq M_h$ holds for all $e^k \in K^k$. We denote by
$||h||_{\infty}$ the supremum of $h$ on the $k$-cells.

The pair $(K,L)$ has {\em bounded geometry }in dimension $n$ if for
all $k\leq n$ there is a bound $M_k$ on the $\ell_1$-norms of the
chains $\partial e^k$, where $e^k$ is a $k$-cell of $K\smallsetminus
L$. Consider a CW-pair $(K,L)$ with bounded geometry in dimension
$n$. It is easy to see that $\delta(C^{k-1}_{(\infty)}(X,Y))\subset
C^k_{(\infty)}(X,Y)$ for $k\leq n$ and hence we can form $H^k_{(\infty)}(X,Y) =
Z^k_{(\infty)}/B^k_{(\infty)}$ if $k\leq n$.

\section{Relative $\ell_{\infty}$-cohomology of $(G,{\mathfrak H})$}

Let $G$ be a group,  $\mathfrak{H} =(H_\lambda)_{\lambda
\in\Lambda}$ a family of subgroups of $G$ and $\mathcal P=\langle X,
{\mathcal H} \mid {\mathcal S},{\mathcal R}\rangle$ a  finite
relative presentation of $G$. We will construct a relative
CW-complex
$K_\mathcal{P}$ associated to the presentation $\mathcal
P$.

For each $\lambda\in\Lambda$ let
$L_\lambda$ be an aspherical CW--complex with one 0-cell
$e^0_\lambda$ such that $\pi_1(L_\lambda,e^0_\lambda)= \widetilde
H_\lambda$.
\begin{remark}
The condition that the $L_\lambda$'s are aspherical is not really
essential. What we need is that for two different choices
$L_\lambda$ and $L'_\lambda$ the isomorphism of the fundamental
groups is induced by a cellular map. So we might choose the
$L_\lambda$'s to be 2-dimensional complexes (see \cite[Lemma
1.5]{ScottWall}).
\end{remark}
We let $L$ denote the disjoint union $L=\sqcup_{\lambda\in\Lambda}
L_\lambda$. Associated to the presentation $\mathcal P$ and $L$
there is a canonical CW-pair $(K_{\mathcal P},L)$ constructed from
$L$ as follows:
\begin{description}
\item[$0$-cells] add one  $0$-cell $e^0$;

\item[$1$-cells] add two types of $1$-cells:
$\{ e^1_\lambda\mid\lambda\in\Lambda\}$ and $\{e^1_x\mid x\in
X\}$.
The $1$-cell $e^1_\lambda$ is attached to $e^0$ and
$e^0_\lambda$ and $e^1_x$ is attached to the $0$-cell $e^0$. We
shall orient the cells $e^1_\lambda$ such that $\partial e^1_\lambda
= e^0_\lambda - e^0$;

\item[$2$-cells] add one $2$-cell $e^2_R$ for each relation $R\in\mathcal R$.
The attaching map for a $2$-cell is given by the corresponding relation.
 \end{description}
 Note that for each $\lambda\in\Lambda$ the subspace
 $L_\lambda\subset K_\mathcal{P}$ is a subcomplex.

By the construction we have $\pi_1(K_\mathcal{P},e^0) = G$ and
$\pi_1(L_\lambda,e^0_\lambda) = \widetilde H_\lambda$. Note also
that
$K_\mathcal{P}^1\cup L\subset K_\mathcal{P}$ is a connected
subcomplex and that $\pi_1(K^1_\mathcal{P}\cup L,e^0)=F$.

Let $K$ be an aspherical CW-complex such that
$K^2=K_\mathcal{P}$.
Such a $K$ can be obtained from $K_\mathcal{P}$
by attaching $k$-cells, $k\geq3$, in order to kill the higher homotopy
groups. Hence we obtain a triple $L\subset K_\mathcal{P} \subset K$.
We consider the universal covering $\pi\co \widetilde{K}\to K$. Note
that $\pi^{-1}(K_\mathcal{P})$ is connected and is hence the
universal covering of $K_\mathcal{P}$. Define
$\overline{L_\lambda}=\pi^{-1}(L_\lambda)$ and
$\overline{L}=\pi^{-1}(L)=\sqcup_{\lambda\in\Lambda}
\overline{L_\lambda}$.

Now, let $(L'_\lambda)_{\lambda\in\Lambda}$ be a second family as
above. Then for each $\lambda\in\Lambda$ there exists a cellular map
$f_\lambda\co L_\lambda\to L'_\lambda$ which induces an isomorphism
between the fundamental groups. Note that we do not require that
$f_\lambda$ is a homotopy equivalence. Starting from the relative
presentation $\mathcal{P}$ we obtain CW-complexes $K'_\mathcal{P}$
and $K'$ by the construction described above. It is obvious that
there exists a cellular map $f\co K\to K'$ extending $f_\lambda$,
$\lambda\in\Lambda$, which is a homeomorphism on each $k$-cell
$e^k\in K\smallsetminus L$, $k=0,1,2$. The lift $\tilde f\co
\widetilde K \to\widetilde{K'}$ induces an isomorphism $\tilde
f_\ast\co C_k(\widetilde K, \overline L)\to C_k(\widetilde{K'},
\overline{L'})$, $k=0,1,2$. More precisely, for  $k=0,1,2$ the image
of each $k$-cell of $\widetilde{K}\smallsetminus \overline{L}$ is a
$k$-cell of $\widetilde{K'}\smallsetminus \overline{L'}$. Hence
$\tilde f$ induces an isomorphism $\tilde f^\ast\co
C_{(\infty)}^k(\widetilde{K'}, \overline{L'}) \to
C_{(\infty)}^k(\widetilde K, \overline L)$ for $k=0,1,2$. By the
symmetry of the construction we obtain in the same way a cellular
map $g\co K'\to K$ and its lift $\tilde g\co
\widetilde{K'}\to\widetilde{K}$ can be chosen such that $(\tilde
f_\ast)^{-1}=\tilde g_\ast\co C_k(\widetilde{K'}, \overline{L'}) \to
C_k(\widetilde K, \overline L)$ for $k=0,1,2$. A direct calculation
gives that $\tilde f^\ast\co H_{(\infty)}^k(\widetilde{K'},
\overline{L'}) \to H_{(\infty)}^k(\widetilde K, \overline L)$ is an
isomorphism for $k=0,1,2$.

\begin{definition}
\label{centrale} Let $G$ be a group which admits a finite
presentation $\langle X, {\mathcal H} \mid {\mathcal S}, {\mathcal
R}\rangle$ relative to a
family of subgroups $\mathfrak H$. The
second $\ell_{\infty}$-cohomology group of $G$ relative to
$\mathfrak H$, denoted by $H^{2}_{(\infty)}(G,{\mathfrak H})$, is
equal to $H^2_{(\infty)}(\widetilde{K},\overline{L})$, where $K$ is
an aspherical CW-complex such that $K^2=K_\mathcal{P}$ and $L\subset
K_\mathcal{P} \subset K$ as above.
\end{definition}

Consider now another finite relative presentation
$\mathcal{P}'=\langle X', {\mathcal H} \mid {\mathcal S},
\mathcal{R}' \rangle$ of $(G,{\mathcal H})$.  The two presentations
are related  by a finite sequence of Tietze transformations (see
\cite{Osin}). As in \cite{Gersten2}, a Tietze transformation is
geometrically realized by an elementary expansion or an elementary
collapse in the Whitehead sens. On the other hand, one can prove
that the group $H^{2}_{(\infty)}(\widetilde K,\overline L)$ does not
change under such an elementary operation (see \cite[Theorem
10.1]{Gersten2}).
\begin{definition}
\label{gersten} Let $G$ be a group which admits a finite
presentation $\langle X, {\mathcal H} \mid {\mathcal S},{\mathcal
R}\rangle$ relative to a family of subgroups $\mathfrak H$. The
second $\ell_{\infty}$-cohomology group of $G$ relative to
$\mathfrak H$ {\em strongly vanishes} if, for some (and hence any)
associated CW-pair $(K,L)$, the sequence $$0 \rightarrow
H^1_{(\infty)}(\widetilde{K},\overline{L}) \rightarrow
C^1_{(\infty)}(\widetilde{K},\overline{L}) /
B^1_{(\infty)}(\widetilde{K},\overline{L})
\overset{\delta}{\rightarrow}
Z^2_{(\infty)}(\widetilde{K},\overline{L}) \rightarrow 0$$ is a
short exact-sequence and there is a bounded section $\sigma \colon
Z^2_{(\infty)}(\widetilde{K},\overline{L}) \rightarrow
C^1_{(\infty)}(\widetilde{K},\overline{L})$ of $\delta$ i.e.\
$\delta\circ\sigma = id$ and there exists $C$ such that
$\|\sigma(z)\|_{\infty}\leq C\|z\|_{\infty}$.
\end{definition}

\subsection{Example}

\label{contrexemple} We give an easy example of a weak relative
hyperbolic group for which the $\ell_{\infty}$-cohomology does not
vanish. Recall that a group $G$ is weakly hyperbolic relative to
$\mathfrak H$ if and only if  the Cayley graph $\Gamma(G,X \cup {\mathcal H})$
is a hyperbolic metric space (see Section~\ref{rh}).

\begin{figure}[htbp]
\resizebox{5cm}{!}{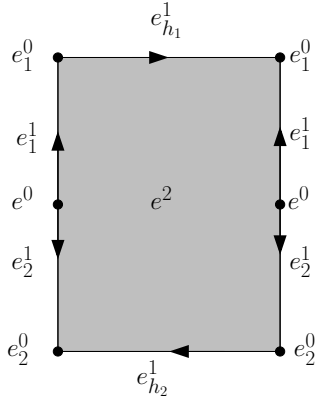}
\caption{The complex $K$.} \protect
\label{fig:pic1}
\end{figure}

Let $G= \mz = \langle h_1,h_2 \mid h_1 h_2 = 1\rangle$ be a finite
presentation of the group $\mz$. Set $\Lambda = \{1,2\}$, $H_i =
\langle h_i\rangle$, $i=1,2$. As in
Section~\ref{rh}, let $\mathcal H$ be the disjoint union of
$\widetilde{H}_1 \setminus \{1\}$ with $\widetilde{H}_2 \setminus
\{1\}$.
Obviously, $\mz$ is weakly hyperbolic
relative to $\mathfrak{H} =(H_\lambda)_{\lambda \in\Lambda}$ since
$\Gamma(G,{\mathcal H})$ is a bounded metric space. The presentation
${\mathcal P} = \langle H_1, H_2 \mid h_1h_2 = 1\rangle$ is a finite
presentation of $\mz$ relative to $\mathfrak H$: $X$ is empty and
$\mathcal R$ is just $h_1h_2=1$. We construct the $2$-complex
$K_\mathcal{P}$ as above. The complex $L_i$ is homeomorphic to the
circle $S^1$ with one $0$-cell $e^0_i$ and one $1$-cell $e^1_{h_i}$.
There is no $1$-cell $e^1_x$ since $X$ is empty so that
$K_{{\mathcal P}}$ is a cylinder: the attaching map of the single
$2$-cell $e^2=e^2_{h_1h_2}$ is given by  the edge path $e^1_1
e^1_{h_1} \overline{e^{1}_1} e^1_2 e^1_{h_2} \overline{e^{1}_2}$.
Here and in the sequel $\overline{e^1_j}$ denotes the $1$-cell
$e^1_j$ with its opposite orientation. The subcomplex $L$ consists
of two disjoint loops and $K$ is aspherical (see
figure~\ref{fig:pic1}).

The universal covering $\pi\co\widetilde{K}\to K$ is homeomorphic
to
the strip $\RR\times [1,2]$ with two boundary componens
$\overline{L_i}= \RR\times \{i\}$, $i=1,2$. We fix a lift of each
cell of $K$ such that $\tilde e^0\in\pi^{-1}(e^0)$, $\tilde
e^0_i\in\pi^{-1}(e^0_i)$, $\tilde e^1_{h_i}\in\pi^{-1}(e^1_{h_i})$,
$i=1,2$, and $\tilde e^2 \in\pi^{-1}(e^2)$ such that
\[\partial \tilde e^1_i = \tilde e^0_i - \tilde e^0,\quad
\partial \tilde e^1_{h_1} =
h_1 \tilde e^0_1 - \tilde e^0_1, \quad
\partial \tilde e^1_{h_2} =
\tilde e^0_2 - h_1 \tilde e^0_2, \]
and
\[\partial \tilde e^2 =
\tilde e^1_1+\tilde e^1_{h_1}-h_1 \tilde e^1_1+
h_1 \tilde e^1_2 + \tilde e^1_{h_2} - \tilde e^1_2\,.
\]

We define a bounded
relative $2$-cocycle $f$ by setting $f(h_1^k\tilde{e}^2) = 1$ for
each $k \in \mz$. Assume now $f = \delta m$ for some relative
$1$-cochain $m$.
The equalities
$ \langle f,h_1^k\tilde{e}^2\rangle = \langle m,h_1^k\partial
\tilde{e}^2\rangle$ and $m(h_1^k \tilde e^1_{h_i})=0$ give
\[1 = m(h_1^k \tilde{e}^1_1)
-m(h_1^{k+1}\tilde{e}^{1}_1) + m(h_1^{k+1}\tilde{e}^{1}_2) -
m(h_1^{k}\tilde{e}^1_2) \text{ for each $k \in \mz$.}\]
By summing from $k=0$ to
$k=n-1$, we get: $n=m(\tilde{e}^1_1) -m(h_1^{n}\tilde{e}^{1}_1) +
m(h_1^{n}\tilde{e}^{1}_2) - m(\tilde{e}^1_2)$.
This implies that
the difference $|m(h_1^n\tilde{e}^{1}_2) -m(h_1^n\tilde{e}^{1}_1)|$
tends toward infinity with $n \to + \infty$. Thus $m$ is not
bounded, so that $f$ is not a bounded coboundary. Therefore
$H^2_{(\infty)}(\mz,{\mathfrak H})$ does not vanish.

\section{Strong vanishing of $H^2_{(\infty)}(G,{\mathfrak H})$ for
strongly relatively hyperbolic groups}

Let $G$ be a group which is strongly hyperbolic relative to a family
of subgroups $ {\mathfrak H} = (H_\lambda)_{\lambda \in \Lambda} $.
As in the previous section we consider a pair $(K,L)$ associated to
a finite presentation $\langle X, {\mathcal H} \mid {\mathcal
    S},{\mathcal R}\rangle$ of $G$
relative to $\mathcal H$. We suppose that $L_{\lambda}$ is the canonical $K(H_{\lambda},1)$
  i.e.\
$L_{\lambda}$ has one $0$-cell, one $1$-cell $e^{1}_{h}$ for each
$h\in H_{\lambda}\setminus\{1\}$
and one $2$-cell $e^{2}_{S}$ for each relation $S\in\mathcal{S}_{\lambda}$.

 Let
$\pi\co\widetilde{K}\to K$ be the universal covering as above. We
equip $\widetilde{K}^{(1)}$ with the following pseudo--metric,
termed the {\em $L$-relative metric}: each $1$-cell of
$\overline{L}$ has length zero, each $1$-cell in
$\pi^{-1}(e^1_\lambda)$, $\lambda\in\Lambda$, has length $1/2$, each
$1$-cell in $\pi^{-1}(e^1_x)$, $x\in X$, has length $1$.
We call {\em cycle} a loop in a graph.

With the notations above:
\begin{definition}
The length of an edge-path $p$ in $\widetilde{K}^{1}$ with respect
to the
$L$-relative metric is termed {\em relative length} and is
denoted by $l_{rel}(p)$.
\end{definition}

Each cycle  $c$ in $\widetilde{K}^{1}$ can be filled by a
singular
disc diagram $D\to \widetilde{K}^2$. The \emph{relative area
$A_{rel}(D)$} of $D$ is the number of $2$-cells in the domain in $D$
which correspond to $2$-cells $e^{2}_{R}$, $R\in\mathcal{R}$. The {\em
relative area $A_{rel}(c)$ of a cycle $c$ in $\widetilde{K}^{1}$} is
the minimal relative area of all diagrams $D$ filling $c$.

An immediate consequence of the definition of strong relative
hyperbolicity is:

\begin{lemma}
\label{trop facile} Suppose that $G$ is strongly hyperbolic relative
to $\mathfrak H$. Then there exists a constant $C \geq 1$ such that,
for any cycle $c$ in $\widetilde{K}^{1}$, $A_{rel}(c) \leq C
l_{rel}(c)$ holds.
\end{lemma}
\begin{proof}
Note that $\Gamma(G,X\cup\mathcal{H})$ is obtained from
$\widetilde{K}^1$ by contracting each edge from
$\pi^{-1}(e_\lambda)$, $\lambda\in\Lambda$.
Let
$\kappa\co\widetilde{K}^1\to\Gamma(G,X\cup\mathcal{H})$ denote the
corresponding surjection.

Now by the very definition we have for each cycle $c$ in $\widetilde{K}^{1}$:
\[ l_\mathit{rel}(c)=l_{X\cup\mathcal{H}}(\kappa(c)) \text{ and }
A_{rel}(c) = \mathit{Area}^{{\mathcal R}}(\kappa(c))\,.\]

The Lemma follows from the existence of a linear relative Dehn-function.
\end{proof}

Let $z \in Z^2_{(\infty)}(\widetilde K, \overline L)$ be a
$\ell_\infty$-relative $2$-cocycle. Let $m \in C^1(\widetilde K,
\overline L)$ be a relative $1$-cochain with $z = \delta m$. Note
that $H^2(\widetilde K, \overline L)=0$ since $H^2(\widetilde K)=0$
and $H^1(\overline L)=0$. From now on, $z$ and $m$ are fixed. We
want to prove the existence of $k \in C^{(1)}_{(\infty)}(\widetilde
K, \overline L)$ such that $z = \delta k$.

Let $P,P'\in\widetilde K^0$. Following \cite{Gersten}, we define a
{\em maximizing path $w(P,P')$} to be a path from $P$ to $P'$
which
maximizes the integer valued function $\nu(\gamma) = \langle
m,\gamma\rangle -  C ||z||_{\infty} l_{rel}(\gamma)$. Here the
maximum is taken among all path $\gamma$ from $P$ to $P'$. The
existence of a maximizing path follows as in \cite[sec.~5]{Gersten}:
\begin{lemma} \cite{Gersten}
\label{copie}
The function $\nu$ always attains its maximum, i.e.\
this maximum is finite and maximizing paths always exist.
\end{lemma}
\begin{proof}
Let $w$, $\gamma$ be two paths in $\widetilde K^1$ from $P$ to
$P'$
and let $D$ be a minimal filling disk of $\gamma^{-1} w$. Then
\[ \langle z, D \rangle = \langle \delta m , D \rangle =
\langle m, w \rangle - \langle m, \gamma \rangle\]
and
\[\langle z, D \rangle \leq \|z\|_\infty A_\mathit{rel}(D)\leq
C\|z\|_\infty \big(l_\mathit{rel}(w) + l_\mathit{rel}(\gamma)\big)
\]
which implies: $\nu(w)= \langle m,w\rangle - C ||z||_{\infty} l_{rel}(w)\leq
\langle m,\gamma\rangle + C \|z\|_{\infty} l_{rel}(\gamma)$.
\end{proof}
In the sequel we fix some cells in $\widetilde K^0$ as follows:
Let $\tilde e^0\in\pi^{-1}(e^0)$,
$\tilde e^0_\lambda\in\pi^{-1}(e^0_\lambda)$,
$\tilde e^1_\lambda\in\pi^{-1}(e^1_\lambda)$, $\lambda\in\Lambda$,
$\tilde e^1_h \in\pi^{-1}(e^1_h)$, $h\in\mathcal H$,
 and
$\tilde e^1_x\in\pi^{-1}(e^1_x)$, $x\in X$, such that
\[\partial \tilde e^1_\lambda =
\tilde e^0_\lambda - \tilde e^0,\quad
\partial \tilde e^1_x =
x \tilde e^0 - \tilde e^0
 \text{ and }
\partial \tilde e^1_h =
h \tilde e^0_\lambda - \tilde e^0_\lambda,\text{ if
$h\in\widetilde
H_\lambda\smallsetminus\{1\}$.}\]

Let $P\in\widetilde K^0$.
We denote by $w(P)$ a  maximizing path  $w(\tilde e^0,P) $
in $\widetilde K^1$ with initial
point $\tilde e^0$ and terminal endpoint $P$. \par From now on we fix for each $\lambda\in\Lambda$ a system of
representatives $\{g^\lambda_i\}_{i \in I_\lambda}$, of $G /
H_\lambda$. For a given $g \in G$, we set
$\overline{g}^\lambda \Def
g^\lambda_i$ if $g H_\lambda = g^\lambda_i H_\lambda$. Hence we have
$\overline{g_1}^\lambda=\overline{g_2}^\lambda$ if and only if $g_1^{-1} g_2\in
H_\lambda$.

\begin{lemma}
\label{ohoh} Let $d \colon C_0(\widetilde{X}) \rightarrow \mz$ be
the $0$-cochain defined by $d(g\tilde{e}^0) = \nu(w(g\tilde{e}^0))$
and $d(g\tilde{e}^0_\lambda) = \nu(w(\overline{g}^\lambda
\tilde{e}^0_\lambda))$. Then $\delta d \in
C^1(\widetilde{K}^1,\overline L)$.
\end{lemma}
\begin{proof}
We check that $\delta d$ vanishes on $C_1(\overline{L})$.
Let $h\in \widetilde H_\lambda\smallsetminus\{1\}$ then
\[(\delta d)(g \tilde{e}^1_h) = d(g\partial \tilde{e}^1_h) =
d(gh \tilde{e}^0_\lambda) - d(g \tilde{e}^0_\lambda) =0\]
because $\overline{gh}^\lambda =
\overline{g}^\lambda$ if $h \in H_\lambda$.
\end{proof}

\begin{proposition}
\label{crucial}
Set $k = -m + \delta d$. Then $k \in
C^1_{(\infty)}(\widetilde{K},\overline{L})$.
\end{proposition}
\begin{proof}
Let $w$ and $w'$ be two edge-paths in $\widetilde K^1$.
If the
terminal vertex of $w$ coincides with the initial vertex of $w'$ we
will simply denote by $ww'$ the composition of the two paths. The
inverse of an edge $g\tilde e^1$ will be denoted by
$\overline{g\tilde e^1}$.

For all $x\in X$ and all $g\in G$ there exist constants $H_1,H_2\geq 0$ such that:
\[
\begin{cases}
\nu(w(g \tilde e^0) g\tilde e^1_x ) +H_1
&= \nu( w(g x\tilde e^0))\\
\nu(w(g x \tilde e^0) \overline{g\tilde e^1_x} ) +H_2
&= \nu( w(g\tilde e^0))\,.
\end{cases}
\]
This is equivalent to:
\[
\begin{cases}
\nu(w(g \tilde e^0)) +\langle m , g\tilde e^1_x\rangle - C\|z\|_\infty +H_1
&= \nu( w(g x\tilde e^0))\\
\nu(w(g x \tilde e^0)) - \langle m, g\tilde e^1_x\rangle - C\|z\|_\infty  +H_2
&= \nu( w(g\tilde e^0))
\end{cases}
\]
and hence $H_1+H_2 = 2C\|z\|_\infty$ and
\begin{align*}
0 &\leq \nu(w(g x\tilde e^0)) -  \nu( w(g \tilde e^0) g\tilde e^1_x ) \\
    &=
\nu(w(g x\tilde e^0)) -  \nu( w(g \tilde e^0))- \langle m , g\tilde
e^1_x\rangle + C\|z\|_\infty \\
    &= H_1\leq H_1+H_2   = 2C\|z\|_\infty\,.\\
\end{align*}
This implies:
\begin{align*} 0\leq k(g\tilde e^1_x) + C\|z\|_\infty &=
\nu(w(g x\tilde e^0)) -  \nu( w(g \tilde e^0)) -\nu( g\tilde e^1_x ) \\
 &= \nu(w(g x\tilde e^0)) -  \nu( w(g \tilde e^0) g\tilde e^1_x )
 \leq 2C\|z\|_\infty
\end{align*}
and therefore $| k(g\tilde e^1_x) | \leq C\|z\|_\infty$.

Note that $w(\overline g^\lambda \tilde e^0_\lambda)$ is an
edge-path from $\tilde e^0$ to $\overline g^\lambda \tilde
e^0_\lambda$ and $g \tilde e^1_\lambda$ is an edge from $g\tilde
e^0$ to $g\tilde e^0_\lambda$. Moreover, $h=  (\overline
g^\lambda)^{-1}g \in H_\lambda$ and that $\overline g^\lambda \tilde
e^1_h$ is an edge from $\overline g^\lambda \tilde e^0_\lambda$ to
$g \tilde e^0_\lambda$. With this notation we obtain that
$w(\overline g^\lambda \tilde e^0_\lambda) (\overline g^\lambda
\tilde e^1_h) (\overline{g \tilde e^1_\lambda})$ is a path from
$\tilde e^0$ to $g\tilde e^0$ and $w(g \tilde e^0)(g\tilde
e^1_\lambda)(\overline{\overline g^\lambda \tilde e^1_h})$ is a path
from $\tilde e^0$ to $\overline g^\lambda \tilde e^0_\lambda$.
Therefore there exist constants $H_3,H_4\geq 0$ such that:
\[
\begin{cases}
\nu(w(g \tilde e^0) )
&= H_3+ \nu( w(\overline g^\lambda \tilde e^0_\lambda)
(\overline g^\lambda \tilde e^1_h) (\overline{g \tilde e^1_\lambda}) )\\
\nu(w( \overline{g}^\lambda \tilde e^0_\lambda ) &= H_4 + \nu(  w(g
\tilde e^0)(g\tilde e^1_\lambda)(\overline{\overline
g^\lambda \tilde
e^1_h})  )\,.
\end{cases}
\]
This is equivalent to:
\[
\begin{cases}
\nu(w(g \tilde e^0) )
&= H_3+ \nu( w(\overline g^\lambda \tilde e^0_\lambda) )
-\langle m , g\tilde e^1_\lambda\rangle - 2C\|z\|_\infty\\
\nu(w( \overline{g}^\lambda \tilde e^0_\lambda ) &= H_4 + \nu(  w(g
\tilde e^0)) + \langle m , g\tilde
e^1_\lambda\rangle
-2C\|z\|_\infty\,.
\end{cases}
\]
As above we abtain
\[ 0\leq \nu(w( \overline{g}^\lambda \tilde e^0_\lambda ) -
\nu( w(\overline g^\lambda \tilde e^0_\lambda) (\overline g^\lambda
\tilde e^1_h) (\overline{g \tilde e^1_\lambda}) )=
H_3\leq
H_3+H_4\leq 4C\|z\|_\infty\] and therefore
\[ k(g\tilde e^1_\lambda) +2C\|Z\|_\infty =
\nu(w( \overline{g}^\lambda \tilde e^0_\lambda ) - \nu( w(\overline
g^\lambda \tilde e^0_\lambda) (\overline g^\lambda \tilde e^1_h)
(\overline{g \tilde e^1_\lambda}) ) \leq 4C\|z\|_\infty \,. \] Here
we have used that $\langle m , g \tilde e^1_h\rangle =0$ for all
$g\in G$ and all
$h\in\mathcal H$.

\end{proof}

\begin{proof}[Proof of Theorem~\ref{sv}]
By definition, $H^2_{(\infty)}(G,{\mathfrak H}) =
H^2_{(\infty)}(\widetilde K,\overline L)$. Let $z$ be a bounded
relative $2$-cocycle of $(\widetilde K ,\overline L)$. Since
$H^2(\widetilde K ,\overline L) = 0$, $z$ is
a relative $2$-coboundary, $z = \delta m$. From Proposition
\ref{crucial}, $k = -m+\delta d$ is a bounded relative $1$-cochain.
But $\delta(-k) = \delta m = z$. Therefore $z$ is a bounded relative
$2$-coboundary. Whence the vanishing of $H^2_{(\infty)}(G,{\mathfrak
H})$. For the strong vanishing, just defining $\sigma(z) = k$
yields the announced section since $\|\sigma(z)\|_\infty = \|k\|_\infty\leq 4C\|z\|_\infty$.
\end{proof}

\section{A converse to the combination theorem}

\label{CCT}

\subsection{Definitions and statement of theorem}

Let $G$ be a group and let $\mathfrak H = {(H_\lambda)}_{\lambda \in
\Lambda}$ be a family of subgroups of $G$.

\begin{convention}
We shall suppose in the sequel that $H_\lambda$ and $H_{\lambda'}$
are not conjugated for $\lambda\neq\lambda'$. Moreover the
$H_\lambda$'s are infinite subgroups.
\end{convention}

\begin{definition}
\label{relative auto} A {\em relative automorphism of $(G,{\mathfrak
H})$} is an automorphism $\alpha$ of $G$ which preserves $\mathfrak
H$ up to conjugacy. More precisely, there is a permutation
$\sigma\in\mathrm{Sym}(\Lambda)$ such that for any $\lambda \in
\Lambda$ there is $g_\lambda \in G$ such that $\alpha(H_\lambda) =
g^{-1}_{\lambda} H_{\sigma(\lambda)} g_{\lambda}$  i.e.\  we have
$i_{g_\alpha} \circ \alpha(H_\lambda) = H_{\sigma(\lambda)}$ where
$i_{g_\alpha}$ is an inner automorphism of $G$. We call $\sigma$ the
permutation \emph{associated} to $\alpha$. If $\sigma$ is the
identity we will say that $\alpha$ {\em fixes $\mathfrak H$ up to
conjugacy}. The group of relative automorphisms will be denoted by
$\mathrm{Aut}(G,\mathfrak{H})$ and the subgroup of relative
automorphisms which fix $\mathfrak H$ up to conjugacy by
$\mathrm{Aut}_0(G,{\mathfrak H})$.
\end{definition}

Let $\mathcal{A} =\{a_1,\ldots,a_n\}$ be a finite set and let
$\F{n}=\langle \mathcal A \rangle$ be the free group with basis
$\mathcal A$. All the free groups considered are finitely generated
free groups. We will denote by $|w|_{\F{n}}$ or $|w|_{{\mathcal A}}$
the word-length of an element of $\F{n}$, depending on whether the
basis $\mathcal A$ has been specified or not. Our convention is that
the distance between $g$ and $h$ in the free group is given by
$|g^{-1}h|_{\F{n}}$.

We suppose that there is an injective homomorphism
$\alpha\co\F{n}\to\mathrm{Aut}(G,\mathfrak{H})$. We will denote
$\alpha_i\Def\alpha(a_i)\in\mathrm{Aut}(G,\mathfrak{H})$ and more
generally  $\alpha_a\Def\alpha(a)$ i.e.\ for all $a,a'\in\F{n}$ we
have $\alpha_{aa'} = \alpha_a\circ \alpha_{a'}$.

We will define a new pair $(G_\mathcal{A},\mathfrak{H}_\mathcal{A})$
in the following
way: we let $G_\mathcal{A}$ denote the semidirect  product
$G_{{\mathcal A}} = G \rtimes_\alpha \F{n}$ i.e.\
$ga\cdot g'a' = g\alpha_a(g')\,aa'$.
For each $\lambda\in\Lambda$ we denote
\[ \mathcal{H}_\lambda \Def \{ ga\in G_\mathcal{A}
\mid \alpha_a(H_\lambda)=g^{-1} H_\lambda g\}\,.\] It is easy to see
that $\mathcal{H}_\lambda\subset G_\mathcal{A}$ is a subgroup
and
that $H_\lambda\subset\mathcal{H}_\lambda$.
\begin{remark}\label{rem:conj}
If there exists $g_b b\in G_\mathcal{A}$ such that
$\alpha_b(H_\lambda)= g_b^{-1} H_{\sigma(\lambda)} g_b$ then the
subgroups $\mathcal{H}_\lambda$
and $\mathcal{H}_{\sigma(\lambda)}$
are conjugate.
\end{remark}
In order to obtain a family of non conjugated subgroups
$\mathfrak{H}_\mathcal{A}$ of $G_\mathcal{A}$ we are proceeding as follows:
\begin{definition}[$\F{n}$-extension of $\mathfrak H$]
\label{extension}
For each $i=1,\ldots,n$, we let $\sigma_i\in\mathrm{Sym}(\Lambda)$
denote the permutation associated to $\alpha_i$ and we let
 $U=\langle
\sigma_1,\ldots,\sigma_n\rangle\subset\mathrm{Sym}$ denote the subgroup
generated by the $\sigma_i$. Let $L\subset\Lambda$ be a system of
orbit representatives i.e.\
\[ \Lambda = \bigsqcup_{\lambda\in L} U\cdot\lambda\,.\]
We now define a  {\em $\F{n}$-extension $\mathfrak{H}_\mathcal{A}$} by $\mathfrak{H}_\mathcal{A}=
(\mathcal{H}_\lambda)_{\lambda\in L}$.
\end{definition}
Note that the groups $\mathcal{H}_\lambda$, $\lambda\in L$, are
uniquely defined up to conjugacy
in the group $G_{{\mathcal A}}$.

Let us now recall the definition of a ``uniform free group of
relatively hyperbolic automorphisms'':
\begin{definition} \cite{GauteroDernier}
Let ${\mathcal P} = \langle X,{\mathcal H} \mbox{ ; } {\mathcal S},
{\mathcal R} \rangle$ be a finite relative presentation of a group
$G$. A {\em uniform free group of relatively hyperbolic
automorphisms of $(G,{\mathfrak H})$} is a free group $\F{n}$
together with a monomorphism $\alpha\co\F{n}\to
\mathrm{Aut}(G,{\mathfrak H})$ for which there exists $\lambda
> 1$, $N,M \geq 1$ such that, for any $g \in G$ with $l_{X \cup {\mathcal H}}(g) \geq M$,
any pair of $a,b \in \F{n}$ with $|a|_{{\F{n}}} = |b|_{{\F{n}}} = N$
and $|a^{-1} b|_{{\F{n}}} = 2 N$ satisfies: $$\lambda l_{X \cup
{\mathcal H}}(g) \leq \mathrm{max}(l_{X \cup {\mathcal
H}}(\alpha_a(g)),l_{X \cup {\mathcal H}}(\alpha_b(g))).$$
\end{definition}

With these definitions in mind, the reader can now go back to
Theorem \ref{un cas particulier}, which is the theorem we are going
to prove. We will however adopt the following:

\begin{convention}
We will assume for the moment that
 $\alpha\co\F{n} \to\mathrm{Aut}_0(G,{\mathfrak H})$.
\end{convention}
As we shall see, the general case is a straightforward implication
of this particular case.

As in Theorem \ref{un cas particulier}, we will assume that $G$ is a group which is
strongly hyperbolic relative to a {\em finite} family ${\mathfrak
H}$ of subgroups ${(H_i)}_{i=1}^k$. Under this finiteness
hypothesis, the automorphisms considered induce quasi isometries on
the group $G$ equipped with the relative metric.

Assume that a finite relative presentation ${\mathcal P} = \langle
X,{\mathcal H} \mbox{ ; } {\mathcal S},{\mathcal R} \rangle$ has
been chosen. Let $(K,L)$ be a CW-pair associated to $\mathcal P$
where each connected component $L_\lambda$ of $L$ is a
$K(H_\lambda,1)$, and such that the $1$-cells in $L$ are in
bijection with the generators of the subgroups in $\mathfrak H$. For
each automorphism $\alpha_i$ we choose a cellular map
$f_i \colon (K,L) \rightarrow (K,L)$ with $(f_i)_{\#} = \alpha_i$ which
fixes the base-point $e^0$ and such that $f_i(L_\lambda) \subset L_\lambda$ for
each connected component $L_\lambda$ of $L$. Let
$K(G_\mathcal{A})$ be the graph
of spaces defined as follows:

\begin{itemize}
  \item the associated combinatorial graph $\Gamma$ is
  the rose with $n$ petals i.e.\ the one point union $\vee_{i=1}^n S^1$  labelled by the $a_i$'s;
  \item the edge and vertex spaces are copies of the complex
  $K$;

Let us recall that, over each open edge of $\Gamma$, $\mathcal G$ is
homeomorphic to $K \times (0,1)$.

  \item the space $K \times \{0\}$ (resp. $K \times \{1\}$)
  associated to the edge with label $a_i$ is glued
  along the vertex space $K$
  by the map $f_i$ (resp.\ by the identity-map).
\end{itemize}

 It is easily checked that $L_\lambda\subset K$ gives rise to a subcomplex
 $L_\lambda(G_{{\mathcal A}})$ and hence a CW-pair, denoted by $(K(G_{{\mathcal
A}}),L(G_{{\mathcal A}}))$. Observe in
 particular that, by construction, each connected component
 $L_\lambda(G_{{\mathcal A}})$ is a $K({\mathcal H}_\lambda,1)$.

Since $G$ is strongly hyperbolic relative to ${\mathfrak H}$, $G$
admits a finite presentation relative to ${\mathfrak H}$. Let
$\langle X, {\mathcal H} \mid {\mathcal R}, {\mathcal S}
\rangle$ be such a finite relative presentation. By definition of ${\mathcal A}$, for each
$a_i \in {\mathcal A}$, for each $H_\lambda \in {\mathfrak H}$,
there is $g_{i,\lambda} \in G$ such that $\alpha_{i}(H_\lambda) =
g^{-1}_{i,\lambda} H_\lambda g_{i,\lambda}$. We denote by
$S_{i,\lambda}$ such a relation. Let ${\mathcal S}^\prime$ be the
union of the relations in $\mathcal S$ with the relations
$S_{i,\lambda}$. Let ${\mathcal R}^\prime$ be the union of the
relations in $\mathcal R$ with the relations $\alpha_{i}(x_j) =
a_i x_j a^{-1}_i$. Then $G_{{\mathcal A}}$ admits $\langle
X,{\mathcal A},{\mathfrak H}_{{\mathcal A}} \mid {\mathcal
R}^\prime, {\mathcal S}^\prime \rangle$ as finite relative
presentation.

\begin{remark}
Constructing a finite relative presentation for $(G_{{\mathcal A}},
{\mathfrak H}_{{\mathcal A}})$ as above is not so hard when
$\alpha(\F{n}) \subset \mathrm{Aut}_{0}(G, {\mathfrak H})$. However,
if the automorphisms $\alpha_i$ only preserve $\mathfrak H$ up to
conjugacy, such a finite relative presentation does not come so
easily without the finite generation of $G$ (or of the
$H_{\lambda}$'s).
\end{remark}

\begin{convention}
In what follows, the CW-pairs $(K,L)$ and $(K(G_{{\mathcal
A}}),L(G_{{\mathcal A}}))$ are graph of spaces as detailed above.
\end{convention}

With this assumption, $({K},{L})$ is canonically embedded in
$({K}(G_{{\mathcal A}}),{L}(G_{{\mathcal A}}))$ since $\Gamma$ has a
unique vertex. We denote by $j \colon (K,L) \rightarrow
({K}(G_{{\mathcal A}}),{L}(G_{{\mathcal A}}))$ this embedding.  As
suggested by the notation, it satisfies $j(L) \subset L(G_{{\mathcal
A}})$ ($j$ is the embedding which induces the canonical injection of
$G$ in $G_{{\mathcal A}}$). The situation is similar for the
universal coverings, which we denote by $\pi \colon \widetilde{K}
\rightarrow K$ and $\pi_{{\mathcal A}} \colon
\widetilde{K}(G_{{\mathcal A}}) \rightarrow K(G_{{\mathcal A}})$
($\pi^{-1}_{{\mathcal A}}(j(K))$ consists of an infinite number of
copies of $\widetilde{K} = \pi^{-1}(K)$).

\begin{definition}
A {\em horizontal edge-path} in $\widetilde{K}(G_{{\mathcal A}})$ is
an edge-path $\gamma$ between two lifts of the base-point ${e^0}$ which
is contained in a connected component $\widetilde K$ of $\pi^{-1}_{{\mathcal
A}}(j(K))$ (the lift, under $\pi_{{\mathcal A}}$, of the complex $K$
canonically embedded in $K(G_{{\mathcal A}})$).

A {\em horizontal geodesic} is a horizontal edge-path which defines
a geodesic of $\widetilde{K}$ equipped with the $\overline{L}$-relative metric.

A {\em corridor ${\mathcal C}_g$} in $\widetilde{K}(G_{{\mathcal
A}})$ is a union of horizontal geodesics which contains, for a given
$g \in G$ and for any $a \in \F{n}$, exactly one horizontal
geodesic, denoted by $\gamma_g(a)$, from $a \tilde e^0$ to $ a
\alpha_{a^{-1}}(g) \tilde e^0$.
\end{definition}

\begin{remark}
A horizontal edge-path in $\widetilde K(G_\mathcal{A})$ projects
under $\pi_\mathcal{A}$ to a closed path representing an element  of
$G \subset G_{{\mathcal A}}$.

A horizontal geodesic is not necessarily (and most often won't be) a
geodesic for $(\widetilde{K}(G_{{\mathcal
A}}),\overline{L}(G_{{\mathcal A}}))$ equipped with the
$\overline{L}(G_{{\mathcal A}})$-relative metric. The fibers
$\widetilde{K}$ in $\pi^{-1}_{{\mathcal A}}(j(K))$ are indeed
distorted (from a geometrical point of view) in the total space.
\end{remark}

\begin{definition}
\label{def intermediaire} A corridor ${\mathcal C}_g$ is {\em
$(\lambda,N,M)$-separated}, with $\lambda
> 1 \mbox{, } N,M \geq 1$, if
for any horizontal geodesic $\gamma_g(w) \in {\mathcal C}_g$ with
$l_{rel}(\gamma_g(w)) \geq M$, any pair of elements $u,v\in \F{n}$
with $|w^{-1}u|_{\F{n}} = |w^{-1}v|_{\F{n}} = N$ and
$|u^{-1}v|_{\F{n}} = 2N$ satisfies:
$$\lambda l_{rel}(\gamma_g(w)) \leq \mathrm{max}(l_{rel}(\gamma_g(wu)),l_{rel}(\gamma_g(wv))).$$
\end{definition}

\begin{remark}
\label{portes ouvertes}
If there exist $\lambda > 1$, $M,N \geq 1$
such that all corridors of $(\widetilde{K}(G_{{\mathcal
A}}),\overline{L}(G_{{\mathcal A}}))$ are $(\lambda, M,
N)$-separated then $\F{n}$ is a uniform free group of relative
automorphisms of $G$.
\end{remark}

The theorem we want to prove is:

\begin{theorem}
\label{general} Let $G$ be a group which is strongly hyperbolic
relative to a finite family $\mathfrak H$ of subgroups. Let
$\alpha\co\F{n}\to \mathrm{Aut}_0(G,{\mathfrak H})$ be a
monomorphism, and let $\mathcal A$ be a basis of $\F{n}$.

If the semi-direct product $G_\mathcal{A}= G \rtimes_\alpha \F{n}$
is strongly hyperbolic relative to ${\mathfrak H}_{\mathcal A}$,
then there exists $\lambda > 1$, $N,M \geq 1$ such that the
corridors of $(\widetilde{K}(G_{{\mathcal
A}}),\overline{L}(G_{{\mathcal A}}))$ are $(\lambda,N,M)$-separated.
\end{theorem}

\begin{remark}
The strong exponential separation property of \cite{GauteroDernier}
involves another condition, which is the exponential separation of
any two vertices representing elements in distinct right $\mathcal
H$-classes, even if the (relative) distance between these vertices
is smaller than the constant $M$. This condition is obviously
necessary, this is most easily seen with Farb's approach
\cite{Farb}: not satisfying this property contradicts the BCP, and
has nothing to do with $\ell_{\infty}$-cohomology. This is why we do
not evoke it in Theorem \ref{general} above.
\end{remark}

\begin{proof}[Proof of Theorem \ref{un cas particulier} assuming Theorem
\ref{general}:] The full statement of Theorem \ref{un cas
particulier}, i.e.\ when $\alpha(\F{n})$ is not necessarily
contained in $\mathrm{Aut}_0(G,{\mathfrak H})$, is deduced from the
following three lemmas. In order to proceed, we fix a monomorphism
$\alpha\co\F{n}\to\mathrm{Aut}(G,{\mathfrak H})$. Since our family
${\mathfrak H}$ is finite it follows that
$\F{0}\Def\alpha^{-1}(\mathrm{Aut}_0(G,{\mathfrak H}))$ is of finite
index in $\F{n}$. We denote by  $\alpha_{0}$ the restriction of
$\alpha$
to $\F 0$. We set ${\mathcal A}$ and ${\mathcal A}_0$ two
basis respectively of $\F{n}$ and $\F 0$.

\begin{lemma}
\label{weeds} With the above natotations: the subgroup $\F 0$ is a finitely
generated subgroup of $\F{n}$ and
its natural embedding in $\F{n}$ defines a quasi isometry between
$\F 0$ and $\F{n}$.
\end{lemma}

\begin{proof} \cite[Proposition 3.19]{GhysHarpe}.
\end{proof}

We define $G_\mathcal{A}= G\rtimes_{\alpha}\F n$ and $G_{\mathcal{A}_0}=
G\rtimes_{{\alpha}_{0}}\F 0$. As in Definition~\ref{extension} let
$\mathfrak{H}_\mathcal{A}$ be the $\F{n}$-extension of
$\mathfrak{H}$ and let $\mathfrak{H}_{\mathcal{A}_0}$ be the
$\F{0}$-extension of $\mathfrak{H}$ i.e.\
\[ \mathfrak{H}_\mathcal{A} = (\mathcal{H}_\lambda)_{\lambda\in L}
\text{ and }
 \mathfrak{H}_{\mathcal{A}_0} =
 (\mathcal{H}^0_\lambda)_{\lambda\in\Lambda}\]
 where $\mathcal{H}^0_\lambda =
 \{ ga\in G_{\mathcal{A}_0} \mid \alpha_a(H_\lambda) =
 g^{-1} H_\lambda g\}$.

\begin{lemma}
\label{michael}
Suppose that $G_\mathcal{A}$ is finitely generated.
If $G_\mathcal{A}$ is strongly hyperbolic relative to
$\mathfrak{H}=(\mathcal{H}_\lambda)_{\lambda\in L}$ then
$G_{\mathcal{A}_0}$ is strongly hyperbolic relative to
$\mathfrak{H}_0=(\mathcal{H}^0_\lambda)_{\lambda\in \Lambda}$.
\end{lemma}
\begin{proof}

Let $\Gamma$ be finitely generated group and suppose that $\Gamma$
is hyperbolic relative to a finite family $(H_j)_{j=1}^k$ of infinite
subgroups.
Moreover,
let $\Gamma_0\subset \Gamma$ be a subgroup of finite index
$p\in\NM$.

We fix a finite system $\{g_{i j}\mid j=1,\ldots, k,\  i=1,\ldots,
p_j\}$
of representatives for the double cosets $G_0 / G \backslash
H_j$ and we define a finite family of subgroups of $\Gamma_0$ by
$H_{ij}\Def(g_{i j} H_j g_{i j}^{-1}) \cap \Gamma_0$. Let
$\mathfrak{H}_0$ be the family $(H_{ij})_{i,j=1}^{ p_j, k}$.

\begin{claim}\label{claim1}
 The group $\Gamma_0$ is hyperbolic relative to  $\mathfrak{H}_0$.
 \end{claim}
\begin{proof}[Proof of the Claim~\ref{claim1}]
Since $\Gamma$ is finitely generated and the family $(H_j)_{j=1}^k$
is finite, each subgroup $H_j$ is finitely generated (see
\cite{Osin} and \cite[Lemma 2.14]{Groves}). Now, the easiest way to
prove the Lemma is to use Bowditch's definition \cite{Bowditch}: a
finitely generated group $\Gamma$ is hyperbolic relative to a family
of finitely generated subgroups $H_1,\ldots, H_k$ if it admits an
action on a hyperbolic graph $K$ such that the following conditions
hold:
\begin{itemize}
\item All edge stabilizers are finite.
\item All vertex stabilizers are finite or conjugate to one of the subgroups $H_1,\ldots, H_k$.
\item The number of orbits of edges is finite.
\item The graph $K$ is fine, that is, for every $n\in \NM$, any edge
of $K$ is contained in finitely many circuits of length n.
\end{itemize}

Now, $\Gamma_0$ is also finitely generated and acts on $K$ and the
edge stabilizers are also finite.
Moreover, there are only finitely
many  $\Gamma_0$-orbits of edges since $\Gamma_0\subset\Gamma$ is of
finite index.

If $\mathrm{Stab}_{\Gamma} (v)$ denotes the stabilizer of the vertex $v$ then
\[ \mathrm{Stab}_{\Gamma_0} (v) =
\mathrm{Stab}_{\Gamma} (v)\cap\Gamma_0\,.\] Note that $
\mathrm{Stab}_{\Gamma_0} (v)$ is infinite if and only if
$\mathrm{Stab}_{\Gamma} (v)$ is infinite since
$\Gamma_0\subset\Gamma$ is of finite index. Moreover, if
$\mathrm{Stab}_{\Gamma_0} (v)$ is infinite there exists $g\in\Gamma$
and $j\in\{1,\ldots,k\}$ such that $g H_j g^{-1} =
\mathrm{Stab}_{\Gamma} (v)$.
Now chose $g_{ij}$ such that $g = g_0
g_{ij} h_j \in G_0 g_{ij} H_j$ hence
\[
g_0 H_{ij} g_0^{-1} = g_0 g_{ij}
H_j g_{ij}^{-1}g_0^{-1} \cap\Gamma_0
= g H_j g^{-1}\cap \Gamma_0 = \mathrm{Stab}_{\Gamma_0} (v)  \,.\]

Hence all vertex stabilizers $ \mathrm{Stab}_{\Gamma_0}(v)$ are finite
or conjugate
to one of the subgroups $H_{ij}$.
Note that if there exists $g_0\in\Gamma_0$ such that
the intersection $g_0 H_{ij} g_0^{-1} \cap H_{i' j'}$ is infinite then
$j=j'$,$i=i'$ and
$g_0\in H_{ij}$.
\end{proof}

The group $G_{\mathcal{A}_0}\subset G_\mathcal{A}$ is of finite
index. Moreover we have for every $\lambda\in L$ that
$\mathcal{H}^0_\lambda = \mathcal{H}_\lambda\cap G_{\mathcal{A}_0}$.
If $\mu\in\Lambda\smallsetminus L$ then there exist $\lambda\in L$
and $ga\in G_\mathcal{A}$ such that $ ga \mathcal{H}_\lambda
a^{-1}g^{-1} = \mathcal{H}_\mu$ (see Remark~\ref{rem:conj}). Hence
\[ ga \mathcal{H}_\lambda a^{-1}g^{-1} \cap G_{\mathcal{A}_0}=
\mathcal{H}_\mu \cap G_{\mathcal{A}_0} = \mathcal{H}_\mu^0\,.\] The
family $\mathfrak{H}_0=(\mathcal{H}^0_\lambda)_{\lambda\in \Lambda}$
is clearly a maximal family of non conjugated subgroups.
\end{proof}

\begin{remark}
\label{moi}
Lemma \ref{michael} is the only place in the whole paper
where we need the finite generation of the semi-direct product
$G_{{\mathcal A}}$. As already evoked in the introduction, we feel
that this could be avoided, but at the expense of an heavy, more
technical work. Let us just sketch the arguments. There are three
steps and the conclusion is given by \cite[Corollary 2.54]{Osin}.
The first step is that $G_{{\mathcal A}_0}$ admits a finite relative
presentation: this is straightforward. The second step is that
$G_{{\mathcal A}_0}$ equipped with its relative metric is
hyperbolic. This is proven by showing that the natural embedding of
$G_{{\mathcal A}_0}$ into $G_{{\mathcal A}}$ is a quasi-isometry.
This is not so hard but here some geometrical arguments are
necessary. They are intuitively clear but writing the details is a
little bit longer. Finally the third and last step is the real
problem: it consists in proving that $G_{{\mathcal A}_0}$ admits a
well-defined relative Dehn function. The idea is to use the fact
that $G_{{\mathcal A}}$ itself admits one, but the real proof
appears difficult to write.
\end{remark}

\begin{lemma}
\label{enfin fini} Assume that $\F 0$ is a uniform free group of
relatively
 hyperbolic automorphisms. Then $\F{n}$ is a uniform free group of
 relatively hyperbolic automorphisms.
\end{lemma}

\begin{proof}

We begin by the following:

\begin{claim}
\label{oneeins}
Assume that $\F{n}$ is
not a uniform free group of relatively hyperbolic automorphisms.
Then for
any sufficiently large
$M > 0$ and $N \geq 1$, for any $\lambda
> 1$, there exist three elements $g,\alpha_{w_0}(g),\alpha_{w_1}(g)$
in $G$ satisfying the following properties: \begin{enumerate} \item
$w_i \in \F 0$, $|w_i|_{\F 0}
\geq N$, $|w^{-1}_1 w_2|_{\F 0} = |w_0|_{\F 0} +
|w_1|_{\F 0}$, \item $l_{rel}(g) \geq M$ and
$l_{rel}(\alpha_{w_i}(g)) < \lambda l_{rel}(g)$.
\end{enumerate}
\end{claim}

\begin{proof}
If $\F{n}$ is not a uniform free group of relatively hyperbolic
automorphisms then:

For any $M
> 0$, $N \geq 1$ and $\lambda
> 1$, there exist three elements $g,\alpha_{w_0}(g),\alpha_{w_1}(g)$
in $G$, with $l_{rel}(g) \geq M$, $|w_i|_{{\mathcal A}} = N$,
$|w^{-1}_1 w_2|_{{\mathcal A}} = 2 N$ and $l_{rel}(\alpha_{w_i}(g))
< \lambda l_{rel}(g)$.

By Lemma \ref{weeds}, $\F 0$ is quasi-isometrically embedded in $\F
n$. Thus for every finite basis ${\mathcal A}_0$ of $\F 0$  there
exists $\mu \geq 1$ such that, for any $w \in  \F 0$:
$$\frac{1}{\mu} |w|_{{\mathcal A}} \leq |w|_{{\mathcal A}_0}.$$

Still by Lemma \ref{weeds}, there is $C > 0$ such that each $w_i$ is
$C$-close, for some positive constant $C$, to an element
$w^\prime_i\in\F 0$. If $N$ is strictly greater than $C$, we can assume
$|{w^\prime_0}^{-1} w^\prime_1|_{{\mathcal A}_0} =
|w^\prime_0|_{{\mathcal A}_0} + |w^\prime_1|_{{\mathcal A}_0}$. \par
From the above two observations,
we get $\mu \geq 1$,
$C \geq 0$ and two elements $w^\prime_i\in\F 0$ with $|{w^\prime_0}^{-1} w^\prime_1|_{{\mathcal A}_0} =
|w^\prime_0|_{{\mathcal A}_0} + |w^\prime_1|_{{\mathcal A}_0}$
such that: $$\frac{1}{\mu} |w|_{{\mathcal A}} - C
\leq |w^\prime_i|_{{\mathcal A}_0}.$$

The automorphisms act by quasi isometry on $G$ equipped with the
$\mathfrak H$-relative metric. Since the distance in $\F n$ from
$w^\prime_i$ to $w_i$ is bounded above by $C$, we get a constant $D
\geq 1$ such that $$l_{rel}(\alpha_{w^\prime_i}(g)) \leq D l_{rel}(\alpha_{w_i}(g)).$$

The proof of the claim readily follows from the above observations.
\end{proof}

As an easy consequence of the definition of a uniform free group of
relatively hyperbolic automorphisms (and using the fact that
automorphisms act by quasi isometries - this is needed for the
existence of $C$ below) we have:

\begin{claim}
\label{twozwei} Assume that $\F 0$
is a uniform free group of relatively hyperbolic automorphisms. Then
there are $M, N \geq 1$, $\lambda > 1$, $C
> 0$ such that, for any $g \in G$ with $l_{rel}(g) \geq M$,
for any integer $j \geq 1$, for any $u,v \in \F 0$ with $|u|_{\F 0}
  \geq jN$, $|v|_{\F 0}
\geq jN$ and $|u^{-1}v|_{\F 0} = |u|_{\F 0} +
|v|_{\F 0}$, $$C \lambda^j l_{rel}(g) \leq
\mathrm{max}(l_{rel}(\alpha_u(g)),l_{rel}(\alpha_v(g))).$$
\end{claim}

Claims \ref{oneeins} and \ref{twozwei} are obviously in
contradiction. We so get the lemma.

\end{proof}

We now complete the proof of Theorem \ref{un cas particulier},
assuming Theorem \ref{general}:

Assume that $G_{{\mathcal A}}$ is strongly hyperbolic relative to
${\mathfrak H}_{{\mathcal A}}$. By Lemma \ref{michael},
$G_{{\mathcal A}_0}$ is strongly hyperbolic relative to ${\mathfrak
H}_{{\mathcal A}_0}$. By Theorem \ref{general} and Remark
\ref{portes ouvertes}, $\F 0 = \langle {\mathcal A}_0 \rangle$ is a
uniform free group of relatively hyperbolic automorphisms. By Lemma
\ref{enfin fini}, $\F{n} = \langle {\mathcal A} \rangle$ is a
uniform free group of relatively hyperbolic automorphisms.
\end{proof}

\begin{remark}
The essential difference between the general case of Theorem \ref{un
cas particulier} and the particular case where the automorphisms fix
each subgroup of $\mathfrak H$ up to conjugacy lies in the fact that
in the first case $G$, equipped with the ``${\mathfrak H}_{{\mathcal
A}} \cap G$-relative metric'', {\em is not} quasi isometrically
embedded in itself by the automorphisms of $\langle {\mathcal A}
\rangle$. This quasi isometric embedding is used in a crucial way
when working with corridors in the next section.
\end{remark}

\section{Proof of Theorem \ref{general}}

We follow the strategy of \cite{Gersten}. In Lemma \ref{evident}
below, we give a key inequality which was proven there for the usual
hyperbolic setting and whose generalization to the relative setting
is straightforward:

\begin{lemma}
\label{evident} Assume that $G$ is strongly hyperbolic relative to
$\mathfrak H$ and let $(K,L)$ be a CW-pair for $(G,{\mathfrak H})$.
Let $z = \delta h \in B^2_{(\infty)}(\widetilde{K},\overline{L})$.
Let $D$ be a filling of a closed edge-path $w$ in $\widetilde{K}^1$.
Then $\langle z,D \rangle \leq ||h||_{\infty} l_{rel}(w)$. In
particular, in the case of strong vanishing of the bounded relative
$2$-cohomology, there exists $C > 0$ such that, for any filling $D$
in $\widetilde{K}^2$, $\langle z,D \rangle \leq C l_{rel}(\partial
D)$.
\end{lemma}

\begin{remark}
\label{je remarque 2} By definition of a corridor $\mathcal {C} \Def
\mathcal {C}_g$ ($g \in G$ is fixed) in $\widetilde{K}(G_{{\mathcal
A}})$, a unique horizontal geodesic $\gamma_g(a)$ in $\mathcal C$ is
associated to each element $a\in\F{n}$. Two geodesics $\gamma \Def
\gamma_g(a)$ and $\gamma' \Def \gamma_g(a')$ of $\mathcal C$ are
called \emph{consecutive} if $|a^{-1} a'|_\mathcal{A} =1$. We
suppose for the moment that $a'=a a_i^\epsilon$. The horizontal path
$\gamma'' = a_i^\epsilon \cdot \gamma$ has the same endpoints as the
horizontal geodesic $\gamma'$. Hence there is a horizontal filling
$D'\Def D^\prime_{\gamma',\gamma''}$ of the loop
$\gamma'(\gamma'')^{-1}$. There is a  loop given by $a_i^\epsilon
\gamma'' a_i^{-\epsilon} (\gamma')^{-1}$. Let  $D''$ denote a
filling of this loop. A \emph{filling} $D_{\gamma,\gamma^\prime}$ of
two consecutive geodesics $\gamma,\gamma^\prime$ in $\mathcal C$ is
defined by concatenating $D'$ and $D''$. By concatenating the
fillings $D_{\gamma,\gamma^\prime}$ so defined for each pair of
consecutive horizontal geodesics in $\mathcal C$, we get a
``filling'' of $\mathcal C$.
\end{remark}

The following lemma is a straightforward generalization, to the
relative setting, of \cite[Proposition 3.1]{Gersten}.

\begin{lemma}
\label{generalisation de Gersten} Let ${\mathcal C} \Def {\mathcal
C}_{g}$ be a corridor in $(\widetilde{K}(G_{{\mathcal
A}}),\overline{L}(G_{{\mathcal A}}))$. Let $u,v$ be two elements of
$\F{n}$ and let $w_1,\cdots,w_n$, $w_1 = u \mbox{, } w_n =v$ be the
elements of $\F{n}$ in the geodesic, in $\Gamma(\F{n},{\mathcal
A})$, from $w_1$ to $w_n$. Let ${\mathcal C}_{u,v}$ be the union of
the horizontal geodesics $\gamma_{i} \Def \gamma_{g}(w_{i})$ in
$\mathcal C$.

There is a filling $D$ of ${\mathcal C}_{u,v}$ and $z \in
Z^2_{(\infty)}(\widetilde{K}(G_{{\mathcal
A}}),\overline{L}(G_{{\mathcal A}}))$ such that
\[
\langle z,D \rangle = \sum^{n-1}_{i=1} l_{rel}(\gamma_i)\,.\]
\end{lemma}

\begin{proof}
We follow the proof of \cite{Gersten}. We first define a bounded relative $1$-cocycle of
$(\widetilde{K},\overline{L})$ by setting, for any horizontal edge
$\tilde e$:
$f(\tilde e) = d_{rel}(\tilde e^0,t(\tilde e)) - d_{rel}(\tilde e^0,i(\tilde e))$,
where $d_{rel}$ denotes the relative distance in
$(\widetilde{K},\overline{L})$. Obviously, when applied to a
geodesic horizontal edge-path $\gamma$ in
$\widetilde{K}(G_{{\mathcal A}})$, we get
$\langle f,\gamma \rangle = l_{rel}(\gamma) = d_{rel}(\tilde
e^0,t(\gamma))$. We define a $2$-cochain $z$ by:

\begin{itemize}
  \item $\langle z,c \rangle = 1$ if the bottom of $c$ is a
  $1$-cell $\tilde{e}^{1}_x$ in ${\mathcal C}_{u,v}$, where
  the $\tilde{e}^{1}_x$ are the lifts, under $\pi_{{\mathcal A}}$, of
  the $1$-cells $e^{1}_{x}$ in $j(K)$ (the image of $K$ in
  $K(G_{{\mathcal A}})$ under its canonical embedding), i.e. the
  $1$-cells associated to the finite set $X$ in the finite relative
  presentation of $G$;
  \item $\langle z,c \rangle = 0$ if the bottom of $c$ is a $1$-cell
  in ${\mathcal C}_{u,v} \cap \pi^{-1}_{{\mathcal A}}(j(L))$, i.e. the
  lift of a
  $1$-cell coming from the relative part of $G$;
  \item $\langle z,c \rangle = \frac{1}{2}$ if the bottom of $c$ is a
  $1$-cell $\tilde{e}^{1}_\lambda$ in ${\mathcal C}_{u,v}$,
  i.e. the lift under $\pi_{{\mathcal A}}$ of a $1$-cell between the
  base point of the complex and the base point of some $L_{\lambda}$;
  \item $\langle z,c \rangle = 0$ if $c$ is any other $2$-cell,
  in particular if $c$ is
  a horizontal $2$-cell.
\end{itemize}

The {\em key-observation} is the following one: if $h$ is a
horizontal edge-path, the value of $z$ on the sum of $2$-cells which
have a $1$-cell of $h$ as bottom is equal to the value of $f$ on
$h$. Since $f$ is a $1$-cocycle of $(\widetilde{K},\overline{L})$,
the product-structure of $K(G_{{\mathcal A}})$ then implies that $z$
vanishes on the $2$-boundaries of $(\widetilde{K}(G_{{\mathcal
A}}),\overline{L}(G_{{\mathcal A}}))$, i.e. $z$ is a $2$-cocycle. By
construction, $z$ vanishes on $\overline{L}(G_{{\mathcal A}})$, i.e.
is a {\em relative} $2$-cocycle. When applying the key-observation
above to the filling $D$ between two consecutive geodesics
$\gamma,\gamma^\prime$ in ${\mathcal C}_{u,v}$ as described in Remark \ref{je
remarque 2}, we get
$\langle z,D \rangle = \langle f,\gamma \rangle =
l_{rel}(\gamma)$. We so get the announced equality.
\end{proof}

Given a corridor ${\mathcal C}_{g}$ and a horizontal geodesic
$\gamma_{g}(w)$ in ${\mathcal C}_{g}$, we say that two horizontal
geodesics $\gamma_{g}(w_{0})$ and $\gamma_{g}(w_{1})$ in ${\mathcal
C}_{g}$ are in a same side of $\gamma_{g}(w)$ if and only if the
geodesic in $\Gamma(\F{n},{\mathcal A})$ from $w_{0}$ to $w_{1}$
does not contain $w$. A {\em side of $\gamma_{g}(w)$} in ${\mathcal
C}_{g}$ is then a maximal union of horizontal geodesics in
${\mathcal C}_{g}$ which are all in a same side of $\gamma_{g}(w)$.

\begin{lemma}
\label{ma methode}
There exists $\lambda_+ > 1$ and $M \geq 1$ such
that, if $\gamma$ is a horizontal geodesic in a corridor $\mathcal C$
with $l_{rel}(\gamma) \geq M$ then there is at least
one side of $\gamma$ in $\mathcal C$ such that any horizontal
geodesic $\gamma^\prime$
in this side satisfies $l_{rel}(\gamma^\prime)
\geq \frac{1}{\lambda_+} l_{rel}(\gamma)$.
\end{lemma}

\begin{proof}
As was already observed, the finiteness of $\mathfrak H$ implies
that the relative automorphisms $\alpha_i$ associated to the $a_i$'s
generating $\F{n}$ act by quasi isometries on $G$ equipped with the
$\mathfrak H$-relative metric. Thus, there is $\mu
> 1$ such that, if $\gamma_0$ is a
horizontal geodesic in $\mathcal C$ consecutive to $\gamma$, then
$l_{rel}(\gamma_0) \geq \frac{1}{\mu} l_{rel}(\gamma)$.

The strong vanishing of
 $H^2_{(\infty)}(G_{{\mathcal A}},{\mathfrak H}_{{\mathcal A}})$
 gives a positive constant $C$ such that $\langle z,D \rangle \leq C l_{rel}(\partial D)$
 ($C$ is the supremum of $||\sigma(z)||_{\infty}$ where $\sigma$
 is the bounded section given by the strong vanishing).

Assume the existence, in $\mathcal C$, of a horizontal geodesic
$\gamma$ such that there exist two horizontal geodesics $\gamma_0
\Def \gamma_{g}(w_{0}),\gamma_1 \Def \gamma_{g}(w_{1})$ in two
distinct sides of $\gamma$ in $\mathcal C$, satisfying the following
properties for some integer $j \geq 1$ (we want to prove that $j$
cannot be chosen arbitrarily large):

\begin{enumerate}
  \item \label{aaaaaa} $l_{rel}(\gamma_i) < \frac{1}{\mu^j} l_{rel}(\gamma)$,
  \item no horizontal geodesic in $\mathcal C$ between $\gamma$ and $\gamma_i$
  satisfies the above inequality,
  \item \label{trois} $\frac{1}{\mu^j} l_{rel}(\gamma) \geq 3 C$.
\end{enumerate}

We consider a filling $D$ of the subset of $\mathcal C$ between
$\gamma_0$ and $\gamma_1$, and a cocycle $z \in
Z^2_{(\infty)}(\widetilde{K}(G_{{\mathcal
A}}),\overline{L}(G_{{\mathcal A}}))$ as given by Lemma
\ref{generalisation de Gersten}. We want to find a minoration of
$\frac{\langle z,D \rangle}{l_{rel}(\partial D)}$ which tends toward
infinity with $j$.

\begin{claim}
$\frac{\langle z,D \rangle}{l_{rel}(\partial
D)}$ is minimal when the $|w_i|_{{\mathcal A}}$'s are minimal.
\end{claim}

\begin{proof}
Let $A_{min}$ and $L_{min}$ be equal to the values respectively of
$\langle z,D \rangle$ and $l_{rel}(\partial D)$ when the
$|w_i|_{{\mathcal A}}$'s are minimal. From (\ref{trois}) above, if
$l_{rel}(\partial D) = L_{min} + x$, then $\langle z,D \rangle \geq A_{min} +
C x$. Thus $\frac{\langle z,D \rangle}{l_{rel}(\partial D)} \geq
\frac{A_{min}}{L_{min}}$ if and only if $C \geq
\frac{A_{min}}{L_{min}}$. As was observed before, this last
assertion is true thanks to the strong vanishing of
$H^2_{(\infty)}(G_{{\mathcal A}},{\mathfrak H}_{{\mathcal A}})$,
which proves the claim.
\end{proof}

\par From our starting observation, Item (\ref{aaaaaa}) implies $|w_i|_{{\mathcal A}} \geq j+1$. By
the claim, we can assume $|w_i| = j+1$. Then $l_{rel}(\partial D)
\leq 4 (j+1) + \frac{2}{\lambda^j} l_{rel}(\gamma)$, whereas $\langle z,D
\rangle \geq (\frac{2}{\lambda} + \cdots + \frac{2}{\lambda^{j-1}} +
\frac{1}{\lambda^j} + 1) l_{rel}(\gamma)$. It follows that the quotient
$\frac{\langle z,D \rangle}{l_{rel}(\partial D)}$ always tends
toward infinity with $j$. As was observed before, this gives a
contradiction with the strong vanishing of the second relative
$\ell_{\infty}$-cohomology group and Lemma \ref{ma methode} follows.
\end{proof}

\begin{proof}[Proof of Theorem \ref{general}]
We argue by contradiction. By a translation in $G_{{\mathcal A}}$,
in order to simplify the notations we can take $w = 1_{\F{n}}$ in
Definition \ref{def intermediaire}. Thus, we assume that, for any
$\lambda
>1$, for any $N,M \geq 1$, there exist $\gamma \Def \gamma_g(1_{\F{n}}),\gamma_0 \Def
\gamma_{g}(w_{0}), \gamma_1 \Def \gamma_{g}(w_{1})$ in a corridor
${\mathcal C} \Def {\mathcal C}_g$ such that $l_{rel}(\gamma) \geq
M$, $l_{rel}(\gamma_i) < \lambda l_{rel}(\gamma)$, $|w_0|_{{\mathcal
A}} = |w_{1}|_{{\mathcal A}} = N$ and $|w^{-1}_0 w_1|_{{\mathcal A}}
= 2N$. From Lemma \ref{ma methode}, we can assume that between
$\gamma_0$ and $\gamma$ all horizontal geodesics of $\mathcal C$
have relative length at least $\frac{1}{\lambda_+} l_{rel}(\gamma)$.
We once again appeal to Lemma \ref{generalisation de Gersten} for
the subset of $\mathcal C$ between $\gamma_0$ and $\gamma$. We get a
bounded relative $2$-cocycle $z$ and a filling $D$ with $\langle z,D
\rangle \geq N \frac{1}{\lambda_+} l_{rel}(\gamma)$. But
$l_{rel}(\partial D) \leq 2 N + (1 + \lambda) l_{rel}(\gamma)$. As
soon as $N,l_{rel}(\gamma)$ are sufficiently large enough and
$\lambda > 1$ sufficiently small enough, we get a contradiction with
the inequality $\langle z,D \rangle \leq C l_{rel}(\partial D)$
given
 by the strong vanishing of
 $H^2_{(\infty)}(G_{{\mathcal A}},{\mathfrak H}_{\mathcal A})$, whence Theorem \ref{general}.
 \end{proof}

\bibliographystyle{plain}
\bibliography{biblio}

\end{document}

%% file: pic1.pdf_t
\begin{picture}(0,0)%
\epsfig{file=pic1.pdf}%
\end{picture}%
\setlength{\unitlength}{3947sp}%
\begingroup\makeatletter\ifx\SetFigFont\undefined%
\gdef\SetFigFont#1#2#3#4#5{%
  \reset@font\fontsize{#1}{#2pt}%
  \fontfamily{#3}\fontseries{#4}\fontshape{#5}%
  \selectfont}%
\fi\endgroup%
\begin{picture}(6157,6345)(376,-5821)
\put(451,-2836){\makebox(0,0)[lb]{\smash{{\SetFigFont{29}{34.8}{\rmdefault}{\mddefault}{\updefault}{\color[rgb]{0,0,0}$e^0$}%
}}}}
\put(451,-3811){\makebox(0,0)[lb]{\smash{{\SetFigFont{29}{34.8}{\rmdefault}{\mddefault}{\updefault}{\color[rgb]{0,0,0}$e^1_2$}%
}}}}
\put(376,-5236){\makebox(0,0)[lb]{\smash{{\SetFigFont{29}{34.8}{\rmdefault}{\mddefault}{\updefault}{\color[rgb]{0,0,0}$e^0_2$}%
}}}}
\put(451,-436){\makebox(0,0)[lb]{\smash{{\SetFigFont{29}{34.8}{\rmdefault}{\mddefault}{\updefault}{\color[rgb]{0,0,0}$e^0_1$}%
}}}}
\put(526,-1786){\makebox(0,0)[lb]{\smash{{\SetFigFont{29}{34.8}{\rmdefault}{\mddefault}{\updefault}{\color[rgb]{0,0,0}$e^1_1$}%
}}}}
\put(4951,-436){\makebox(0,0)[lb]{\smash{{\SetFigFont{29}{34.8}{\rmdefault}{\mddefault}{\updefault}{\color[rgb]{0,0,0}$e^0_1$}%
}}}}
\put(5101,-2836){\makebox(0,0)[b]{\smash{{\SetFigFont{29}{34.8}{\rmdefault}{\mddefault}{\updefault}{\color[rgb]{0,0,0}$e^0$}%
}}}}
\put(4951,-3811){\makebox(0,0)[lb]{\smash{{\SetFigFont{29}{34.8}{\rmdefault}{\mddefault}{\updefault}{\color[rgb]{0,0,0}$e^1_2$}%
}}}}
\put(5026,-5236){\makebox(0,0)[lb]{\smash{{\SetFigFont{29}{34.8}{\rmdefault}{\mddefault}{\updefault}{\color[rgb]{0,0,0}$e^0_2$}%
}}}}
\put(4951,-1711){\makebox(0,0)[lb]{\smash{{\SetFigFont{29}{34.8}{\rmdefault}{\mddefault}{\updefault}{\color[rgb]{0,0,0}$e^1_1$}%
}}}}
\put(2476,-5686){\makebox(0,0)[lb]{\smash{{\SetFigFont{29}{34.8}{\rmdefault}{\mddefault}{\updefault}{\color[rgb]{0,0,0}$e^1_{h_2}$}%
}}}}
\put(2701,164){\makebox(0,0)[lb]{\smash{{\SetFigFont{29}{34.8}{\rmdefault}{\mddefault}{\updefault}{\color[rgb]{0,0,0}$e^1_{h_1}$}%
}}}}
\put(2701,-2836){\makebox(0,0)[lb]{\smash{{\SetFigFont{29}{34.8}{\rmdefault}{\mddefault}{\updefault}{\color[rgb]{0,0,0}$e^2$}%
}}}}
\end{picture}%